\documentclass[a4paper,11pt,leqno]{amsart}
\usepackage[ansinew]{inputenc}
\usepackage{graphicx}
\usepackage{picins}

\usepackage{url}

\usepackage[colorlinks=true,linkcolor=blue,citecolor=magenta,urlcolor=red]{hyperref}

\hypersetup{pdfinfo={
  Title=Transcriptions of Six Letters of Bernhard Riemann Preserved in the SIL,
  Author=Wolfgang Gabcke,
  Subject=Transcriptions of six letters of Bernhard Riemann preserved in the SIL,
  Keywords={Bernhard Riemann, Elise Riemann, Wilhelm Weber, letters, German transcription, Smithsonian Libraries}},
            linktocpage
}

\newcommand\mycleardoublepage{%
   \clearpage\ifodd \value{page} \else \thispagestyle{empty}\null\newpage\fi
}

\theoremstyle{plain}

\theoremstyle{definition}

\theoremstyle{remark}

\usepackage{scrextend}

\begin{document}
%
\title[TRANSCRIPTIONS OF SIX LETTERS OF BERNHARD RIEMANN PRESERVED IN THE SIL]
      {TRANSCRIPTIONS OF SIX LETTERS\\OF BERNHARD RIEMANN PRESERVED IN THE SIL}
\author[WOLFGANG GABCKE]{WOLFGANG GABCKE}

\date{26\textsuperscript{th} February 2016; 1\textsuperscript{st} edition}

\begin{abstract}
The German transcriptions of six letters of Bernhard Riemann preserved in the Smithsonian Libraries (SIL) in
Washington DC are given.
\end{abstract}
%
%
\maketitle
%
\deffootnotemark{\textsuperscript{(\thefootnotemark)}}
\deffootnote{2em}{1.6em}{(\thefootnotemark)\enskip}
%
\tableofcontents
\vspace{-4mm}
%
%
%
%
%
%
\section*{Introduction}
\noindent
There is a collection of six letters of Bernhard Riemann preserved in the Smithsonian Libraries (SIL) in Washington
DC\footnote{See \cite{Smi} and for the supplement of the collection, figures \ref{Abb01} and \ref{Abb02}.},
a donation of the Burndy Library DSI in 1974\footnote{In this year, a quarter of the Burndy Library
collection came to the SIL.}. It is possible that there are other manuscripts of Bernhard Riemann
in the SIL or in the Huntington Library in San Marino Ca., where the remaining collection of the Burndy Library is
today. Therefore, it would be very important to find out how this collection has come from Germany to the Burndy
Library before 1974.

With the exception of the first letter, which was written in 1857 in Bremen, the other letters were written in
Göttingen between 1857 and 1866. The facsimiles clearly show that the letters were written on expensive silk paper
except of the second one, but all of them in neat handwriting. Obviously, these letters are not drafts but the
originals that have been sent to the recipients. According to a private communication of Erwin Neuenschwander, it
is very likely that all letters were sent to the same person\footnote{Riemann never uses a name in the
welcome clause of these letters.}, namely Wilhelm Weber\footnote{Wilhelm Eduard Weber (1804--1891), German
physicist.}, a fatherly friend of Bernhard Riemann. This still needs to be verified.

Unlike the first five letters, the sixth letter\footnote{See figure \ref{Abb03} for its first page.}of April 1866
is not in Riemann's hand. It may have been dictated
by Riemann and written by his wife Elise\footnote{Elise Riemann, born Koch (1835--1904).} perhaps because he was
too weak to do it himself.\footnote{Riemann died three months later of his tuberculosis disease.} Only the last
line with date and signature is in his hand. Therefore, I made a comparison of the handwriting of this letter with
that of an autograph of Elise Riemann, a letter to Richard Dedekind\footnote{Richard Dedekind (1831--1916), German
mathematician.}, written 27 November 1871 (see \cite{Eli} and figures \ref{Abb04}, \ref{Abb05}) and
located at the \textit{Nie\-der\-säch\-si\-sche Staats- und Uni\-ver\-si\-täts\-bi\-blio\-thek Göt\-tin\-gen}
(SUB). From this investigation it is evident that Elise Riemann indeed wrote the sixth letter since the
handwriting in both letters is identical.

Layout, spelling and punctuation of the handwritten originals are maintained in the German transcriptions,
regardless of obvious errors.
%
%
%
%
%
%
\newpage
\section{Letter of 7\texorpdfstring{\textsuperscript{th}}{th} July 1857}
\noindent
\large
\phantom{aa}\\[1.0ex]
\phantom{aa}\\[1.0ex]
\phantom{aaaaaaaaaaaaaaaaaaaa}Hochgeehrtester Herr Professor!\\[1.0ex]
\phantom{aa}\\[1.0ex]
\phantom{aa}\\[1.0ex]
\phantom{aaaa}Nachdem ich meine Anhandlung über die Abel'schen Funktionen hier\\[1.0ex]
\phantom{aa}glücklich vollendet hatte, beabsichtigte ich heute nach Göttingen zurück\\[1.0ex]
\phantom{aa}zu reisen und würde Sie dann mündlich in so manchen Dingen, wo\\[1.0ex]
\phantom{aa}ich Ihres Rathes bedarf, um Verhaltungsregeln gebeten haben. Ich erhalte\\[1.0ex]
\phantom{aa}aber nun eben einen Brief von Hrn. Dr. Borchardt\footnote{Karl Wilhelm Borchardt (1817--1880),
            German mathematician.}, worin dieser\\[1.0ex]
\phantom{aa}über einige Stellen in einer schon im Druck begriffenen Abhandlung\\[1.0ex]
\phantom{aa}umgehend Auskunft erbittet und der mich daher zwingt, meine Ab-\\[1.0ex]
\phantom{aa}reise noch ein paar Tage hinauszuschieben.\\[1.0ex]
\phantom{aaaa}Nach einem Briefe des Hrn. Dr. Dedekind\footnote{Richard Dedekind (1831--1916),
            German mathematician.} ist es nun\(\!\!\!\!\!\!\!\!\!\!\)\raisebox{-0.6ex}{\huge\textbf{---}}
            die höchste Zeit,\\[1.0ex]
\phantom{aa}die Vorlesungen für das nächste Semester anzukündigen. Da ich nun\\[1.0ex]
\phantom{aa}------------------------------------------ {\small{new page}} --------------------------------------------\\[1.0ex]
\phantom{aa}ohne Ihren Rath mich über meine Vorlesungen nicht zu entscheiden wage,\\[1.0ex]
\phantom{aa}weil ich nicht weiß, wie viel Zeit ich im nächsten Semester nach Ihrer Ab-\\[1.0ex]
\phantom{aa}sicht auf Vorlesungen verwenden darf, und ob Sie mehr zu einem Gegen-\\[1.0ex]
\phantom{aa}stande der mathematischen Physik oder einem rein mathematischen ra-\\[1.0ex]
\phantom{aa}then, so hoffe ich, Sie werden mir gütigst verzeihen, wenn ich Sie bitte,\\[1.0ex]
\phantom{aa}nöthigenfalls für mich eine Entscheidung zu treffen und Ihnen zu diesem\\[1.0ex]
\phantom{aa}Zwecke meine Pläne mittheile. Für den Fall, daß ich nur ein einstündiges\\[1.0ex]
\phantom{aa}Publicum lesen soll, welches nur wenig Zeit zur Vorbereitung erfordert,\\[1.0ex]
\phantom{aa}würde ich etwa problemata physica selecta ankündigen. Rathen Sie\\[1.0ex]
\phantom{aa}zu einer etwas größeren Vorlesung aus der mathematischen Physik,\\[1.0ex]
\phantom{aa}so würde ich die Theorie der Elasticität fester Körper nach Lamé\footnote{Gabriel Lamé (1795--1870),
            French mathematician and physicist.},\\[1.0ex]
\phantom{aa}oder mathematische Optik, falls diese nicht schon gelesen wird, wählen\\[1.0ex]
\phantom{aa}oder auch, um mich im Vortrage von Gegenständen der mathe-\\[1.0ex]
\phantom{aa}matischen Physik zu üben, einen elementareren Theil derselben be-\\[1.0ex]
\phantom{aa}handeln. Sehr gern würde ich auch eine meiner Vorlesungen aus\\[1.0ex]
\phantom{aa}den letzten Semestern (über elliptische und Abel'sche Functionen und\\[1.0ex]
\phantom{aa}über hypergeometrische Reihen) zum zweiten Male lesen, um mehr\\[1.0ex]
\phantom{aa}Sorgfalt auf den Vortrag verwenden zu können.\\[1.0ex]
\phantom{aa}------------------------------------------ {\small{new page}} --------------------------------------------\\[1.0ex]
\phantom{aa}Ich bitte um die herzlichsten Grüße an die lieben Ihren,\\[1.0ex]
\phantom{aaaa}Ihr in dankbarster Hochachtung\\[1.0ex]
\phantom{aaaaaaaaaaaaaaaaaaaaaaaaaaaaaaaaa}Ihnen\\[1.0ex]
\phantom{aaaaaa}Bremen 7.\ Juli\phantom{aaaaaaaaaaaaaaaaaaaaaaaa}treu ergebener\\[1.0ex]
\phantom{aaaaaaaaa}1857\phantom{aaaaaaaaaaaaaaaaaaaaaaaaaaaaaii}B.\ Riemann\\[1.0ex]
\normalsize
%
%
%
%
%
\newpage
\section{Letter of 11\texorpdfstring{\textsuperscript{th}}{th} July 1857}
\noindent
\large
\phantom{aaaa}\\[1.0ex]
\phantom{aaaa}\\[1.0ex]
\phantom{aaaaaaaaaaaaaaaaaaa}Hochgeehrtester Herr Professor\\[1.0ex]
\phantom{aaaa}\\[1.0ex]
\phantom{aaaa}\\[1.0ex]
\phantom{aaaa}Nach Empfang Ihres so überaus gütigen Briefes vom 9.\ d.\ M.\footnote{Short for ``des Monats''.}\\[1.0ex]
\phantom{aaaa}bin ich sofort hieher geeilt, und ich würde mich Ihnen und den ver-\\[1.0ex]
\phantom{aaaa}ehrten Ihrigen schon präsentiert haben, wenn ich nicht gewünscht hätte,\\[1.0ex]
\phantom{aaaa}den mir freundlichst übersandten Brief Ihrer verehrten Frau\\[1.0ex]
\phantom{aaaa}Schwägerin\footnote{Wife of Ernst Heinrich Weber (1795--1878) or Eduard Friedrich Weber (1806--1871),
              German physiologists and anatomists.} vorher zu beantworten.\\[1.0ex]
\phantom{aaaaaaa}Mit der Bitte, Ihr diese Antwort gütigst zu übergeben\\[1.0ex]
\phantom{aaaa}\\[1.0ex]
\phantom{aaaaaaaaaaaaaaaaaaaaaaaaaaaa}Ihr\\[1.0ex]
\phantom{aaaa}\\[1.0ex]
\phantom{aaaaaaaaaaaaaaaaaaaaaaaaaaaaaaaa}in dankbarster Hochachtung Ihnen\\[1.0ex]
\phantom{aaaa}Göttingen, den 11.\ Juli\phantom{aaaaaaaaaaaaaaaaaa}ergebener\\[1.0ex]
\phantom{aaaaaaaaaaa}1857.\phantom{aaaaaaaaaaaaaaaaaaaaaaaaaaa}Riemann\\[1.0ex]
\normalsize
%
%
%
%
%
%
\newpage
\section{Letter of 8\texorpdfstring{\textsuperscript{th}}{th} November 1861}
\noindent
\large
\phantom{aaaaaaaaa}\\[1.0ex]
\phantom{aaaaaaaaa}\\[1.0ex]
\phantom{aaaaaaaaaaaaaaaaaaaa}Hochgeehrter Herr Professor!\\[1.0ex]
\phantom{aaaaaaaaa}\\[1.0ex]
\phantom{aaaaaaaaa}Mit der Bitte, die Verzögerung gütigst zu entschuldigen\\[1.0ex]
\phantom{aaaaaaaaa}übersende ich Ihnen hieneben die Vorschläge für die Wah-\\[1.0ex]
\phantom{aaaaaaaaa}len der Societät\footnote{Königliche Gesellschaft der Wissenschaften zu Göttingen.}.
                   Die gemeinschaftliche Wahl der beiden\\[1.0ex]
\phantom{aaaaaaaaa}von mir vorgeschlagenen Herren zu Correspondenten\\[1.0ex]
\phantom{aaaaaaaaa}dürfte sich besonders durch die enge Beziehung empfehlen,\\[1.0ex]
\phantom{aaaaaaaaa}in der ihre Arbeiten zu einander stehen. Den bespro-\\[1.0ex]
\phantom{aaaaaaaaa}chenen Vorschlag des Herrn Brioschi\footnote{Francesco Brioschi (1824--1897),
            Italian mathematician.} habe ich für dies-\\[1.0ex]
\phantom{aaaaaaaaa}mal unterlassen, um die Anzahl der Wahlvorschläge\\[1.0ex]
\phantom{aaaaaaaaa}nicht noch mehr zu vergrößern.\\[1.0ex]
\phantom{aaaaaaaaaaaa}Ganz\\[1.0ex]
\phantom{aaaaaaaaa}\\[1.0ex]
\phantom{aaaaaaaaa}\\[1.0ex]
\phantom{aaaaaaaaa}Göttingen, den 8.\ November\phantom{aaaaaaaaaaa}der Ihrige\\[1.0ex]
\phantom{aaaaaaaaaaaaaaaa}1861.\phantom{aaaaaaaaaaaaaaaaaaaaaaaai}B.\ Riemann\\[1.0ex]
\normalsize
%
%
%
%
%
%
\newpage
\section{Letter of 29\texorpdfstring{\textsuperscript{th}}{th} June 1863}
\noindent
\large
\phantom{aaaaaa}\\[1.0ex]
\phantom{aaaaaa}\\[1.0ex]
\phantom{aaaaaaaaaaaaaaaaaaaa}Hochzuverehrender Herr Hofrath!\\[1.0ex]
\phantom{aaaaaa}\\[1.0ex]
\phantom{aaaaaa}\\[1.0ex]
\phantom{aaaaaaaaa}Soeben erhielt ich von Hr. Prof. Schering\footnote{Ernst Christian Julius Schering (1833--1897),
                   German mathematician.} einige Manuscripte\\[1.0ex]
\phantom{aaaaaa}zu G.\footnote{Carl Friedrich Gauß (1777--1855), German mathematician and astronomer.}
                Werken zugesandt, mit den beifolgenden Zeilen begleitet.\\[1.0ex]
\phantom{aaaaaa}Nach Ihren Äußerungen über diese Angelegenheit glaubte ich bei\\[1.0ex]
\phantom{aaaaaa}Ihnen anfragen zu müssen, was ich hierauf thun soll.\\[1.0ex]
\phantom{aaaaaaaaa}Da uns Allen an der Beschleunigung des Druckes gelegen ist, die\\[1.0ex]
\phantom{aaaaaa}von Hrn. Prof. Schering mir gestellte Frist aber zu einer Durchsicht\\[1.0ex]
\phantom{aaaaaa}des Manuscripts offenbar ganz unzureichend ist, so würde ich ihm\\[1.0ex]
\phantom{aaaaaa}das Manuscript morgen früh ungelesen wieder zugesandt haben,\\[1.0ex]
\phantom{aaaaaa}wenn ich nicht nach Ihrer Mittheilung von neulich hätte glauben\\[1.0ex]
\phantom{aaaaaa}müssen, daß eine nochmalige Revision wünschenswerth sei.\\[1.0ex]
\phantom{aaaaaa}----------------------------------------- {\small{new page}} ---------------------------------------\\[1.0ex]
\phantom{aaaaaa}Ich erlaube mir nun Sie zu bitten, mich mit ein paar Zeilen\\[1.0ex]
\phantom{aaaaaa}zu benachrichtigen, ob Sie eine nochmalige Durchsicht trotz der\\[1.0ex]
\phantom{aaaaaa}jetzt dabei unvermeidlichen Verzögerung des Drucks noch wünschen\\[1.0ex]
\phantom{aaaaaa}So gern ich diese im Interesse der Sache übernehme, so möchte ich\\[1.0ex]
\phantom{aaaaaa}doch um keinen Preis, wenn es nicht nöthig ist, hemmend\\[1.0ex]
\phantom{aaaaaa}eingreifen.\\[1.0ex]
\phantom{aaaaaaaaa}In dankbarer Ergebenheit\\[1.0ex]
\phantom{aaaaaa}\\[1.0ex]
\phantom{aaaaaaaaaaaaaaaaaaaaaaaaaaaaaaaaaaaa}Ihr\\[1.0ex]
\phantom{aaaaaa}Göttingen, den 29.\ Juni\phantom{aaaaaaaaaaaaaaaaaaaaaa}B.\ Riemann\\[1.0ex]
\phantom{aaaaaaaaaaaaa}1863.\\[1.0ex]
\normalsize
%
%
%
%
%
%
\newpage
\section{Letter of 29\texorpdfstring{\textsuperscript{th}}{th} July 1863}
\noindent
\large
\phantom{aaa}\\[1.0ex]
\phantom{aaa}\\[1.0ex]
\phantom{aaaaaaaaaaaaaaaaaaaaaa}Verehrtester Herr Hofrath!\\[1.0ex]
\phantom{aaa}\\[1.0ex]
\phantom{aaa}Einliegend sende ich Ihnen den Brief von Dedekind\footnote{Richard Dedekind (1831--1916),
             German mathematician.} wieder zu, für\\[1.0ex]
\phantom{aaa}dessen gütige Mittheilung dankend; ich bitte zugleich um Entschuldigung,\\[1.0ex]
\phantom{aaa}wenn ich Ihnen die Manuscripte noch nicht mit der Erklärung, daß\\[1.0ex]
\phantom{aaa}sie druckfertig seien, wieder zustellen kann. Dedekind hatte die von ihm\\[1.0ex]
\phantom{aaa}gewünschten Abänderungen noch nicht in das für den Setzer bestimmte\\[1.0ex]
\phantom{aaa}Manuscript aufgenommen, sondern sie als von ihm vorgeschlagen\\[1.0ex]
\phantom{aaa}auf einem besondern Blatte zusammengestellt. Schering\footnote{Ernst Christian Julius Schering
             (1833--1897), German mathematician.}, der heute\\[1.0ex]
\phantom{aaa}Morgen bei mir war, hat mich nun gebeten, ihm das Manuscript zur\\[1.0ex]
\phantom{aaa}Durchsicht zu übersenden, um sich dann mit mir, bevor es bei der\\[1.0ex]
\phantom{aaa}Commission ciruliert, besprechen zu können; und ich werde dies\\[1.0ex]
\phantom{aaa}thun, nachdem ich es selbst gelesen habe, in der Hoffnung daß Sie\\[1.0ex]
\phantom{aaa}damit einverstanden sind, oder mich gütigst benachrichtigen, wenn\\[1.0ex]
\phantom{aaa}Sie es anders wünschen sollten.\\[1.0ex]
\phantom{aaa}-------------------------------------------- {\small{new page}} -----------------------------------------\\[1.0ex]
\phantom{aaa}Ich bedaure sehr, daß hiedurch eine neue Verzögerung des Drucks\\[1.0ex]
\phantom{aaa}um ein paar Tage eintritt, die Ihnen gewiß nicht lieb sein wird;\\[1.0ex]
\phantom{aaa}doch lag es nicht in meiner Macht sie zu vermeiden, und ich werde\\[1.0ex]
\phantom{aaa}Alles thun, sie möglichst kurz zu machen.\\[1.0ex]
\phantom{aaaaaa}In dankbarer Ergebenheit\\[1.0ex]
\phantom{aaa}\\[1.0ex]
\phantom{aaaaaaaaaaaaaaaaaaaaaaaaaaaaaaaaaa}Ihr\\[1.0ex]
\phantom{aaa}Göttingen den 29.\ Juli\phantom{aaaaaaaaaaaaaaaaaaaaaa}Riemann\\[1.0ex]
\phantom{aaaaaaaaa}1863.\\[1.0ex]
\normalsize
%
%
%
%
%
%
\newpage
\section{Letter of 18\texorpdfstring{\textsuperscript{th}}{th} April 1866}
\noindent
\large
\phantom{aaaaaaaaa}\\[1.0ex]
\phantom{aaaaaaaaa}\\[1.0ex]
\phantom{aaaaaaaaaaaaaaaaa}Hochgeehrter Herr Geheimer Hofrath.\\[1.0ex]
\phantom{aaaaaaaaa}\\[1.0ex]
\phantom{aaaaaaaaa}Sie haben sich, wie in so vielen Dingen, so auch bei mei\\[1.0ex]
\phantom{aaaaaaaaa}ner Untersuchung über das Ohr mit so großer Güte meiner\\[1.0ex]
\phantom{aaaaaaaaa}angenommen, daß ich es wage, Sie, so wenig ich mich auch\\[1.0ex]
\phantom{aaaaaaaaa}durch die Sache selbst dazu berechtigt halte, mit einem\\[1.0ex]
\phantom{aaaaaaaaa}Briefe darüber zu beschäftigen.\\[1.0ex]
\phantom{aaaaaaaaa}Ihr gütiges Anerbieten mir anatomische Präparate\\[1.0ex]
\phantom{aaaaaaaaa}von Herrn Hofrath Haenle\footnote{Jakob Henle (1809--1885), German anatomist, pathologist and
                   physician; Mrs. Riemann uses the wrong spelling ``Haenle''.} zu verschaffen ist mir zwar\\[1.0ex]
\phantom{aaaaaaaaa}sehr angenehm, weil ich dadurch eine deutlichere An{\raisebox{1.5pt}{\tiny=}}\\[1.0ex]
\phantom{aaaaaaaaa}schauung von allen das Ohr betreffenden Verhältnissen\\[1.0ex]
\phantom{aaaaaaaaa}gewinnen kann. Aber direct kann es zur Förderung mei{\raisebox{1.5pt}{\tiny=}}\\[1.0ex]
\phantom{aaaaaaaaa}ner Arbeit nichts beitragen. Wenn ich mir eine\\[1.0ex]
\phantom{aaaaaaaaa}Gunst von Herrn Hofrath Haenle erbitten möchte, so wäre\\[1.0ex]
\phantom{aaaaaaaaa}es die, daß er meine Abhandlung vor dem Druck durch\\[1.0ex]
\phantom{aaaaaaaaa}läse und mir die Stellen bezeichnete, wo ich etwa\\[1.0ex]
\phantom{aaaaaaaaa}eine nach seinen Beobachtungen unrichtige oder\\[1.0ex]
\phantom{aaaaaaaaa}unwahrscheinliche Hypothese gemacht hätte. Ich glaube\\[1.0ex]
\phantom{aaaaaaaaa}zwar die Bedingungen, denen der Apparat und\\[1.0ex]
\phantom{aaaaaaaaa}seine Theile genügen müssen, vollständiger er{\raisebox{1.5pt}{\tiny=}}\\[1.0ex]
\phantom{aaaaaaaaa}----------------------------------- {\small{new page}} --------------------------------------\\[1.0ex]
\phantom{aaaaaaaaa}kannt und nach der Erfahrung bestimmt zu haben,\\[1.0ex]
\phantom{aaaaaaaaa}als dies bisher gelungen war aber bei dem Aufsuchen\\[1.0ex]
\phantom{aaaaaaaaa}der einzelnen Vorrichtungen, durch welche diese Bedin{\raisebox{1.5pt}{\tiny=}}\\[1.0ex]
\phantom{aaaaaaaaa}gungen erfüllt sind, kann man außerordentlich leicht\\[1.0ex]
\phantom{aaaaaaaaa}fehlgreifen, wenn die anatomischen Angaben für diesen\\[1.0ex]
\phantom{aaaaaaaaa}Zweck nicht genau genug sind, und besonders wenn\\[1.0ex]
\phantom{aaaaaaaaa}man sie so wenig vollständig kennt, wie ich.\\[1.0ex]
\phantom{aaaaaaaaa}Besonders schwierig ist es mir über die Form des\\[1.0ex]
\phantom{aaaaaaaaa}Steigbügels und die Befestigung des Hammers im\\[1.0ex]
\phantom{aaaaaaaaa}Paukenfelle die nöthigen Data mir zu verschaffen.\\[1.0ex]
\phantom{aaaaaaaaa}Ueber die Befestigung des Hammers können ältere\\[1.0ex]
\phantom{aaaaaaaaa}Präparate gar keinen Aufschluß geben, sondern nur\\[1.0ex]
\phantom{aaaaaaaaa}solche von ganz frisch geschlachteten Thieren. Es ist mir\\[1.0ex]
\phantom{aaaaaaaaa}ja leicht, mir Köpfe von ganz frisch geschlachtetem\\[1.0ex]
\phantom{aaaaaaaaa}Vieh holen zu lassen; und ich habe mir auf diese\\[1.0ex]
\phantom{aaaaaaaaa}Weise bei dem Steigbügel Licht verschafft. Die\\[1.0ex]
\phantom{aaaaaaaaa}Befestigung des Hammers wurde aber durch die rohen\\[1.0ex]
\phantom{aaaaaaaaa}Instrumente die mir zum Zerbrechen der das Ohr\\[1.0ex]
\phantom{aaaaaaaaa}einschließenden Knochen zu Gebote standen, immer\\[1.0ex]
\phantom{aaaaaaaaa}zerstört.\\[1.0ex]
\phantom{aaaaaaaaa}Sollte Ihr Herr Neffe\footnote{What is his name?} mich nun noch besuchen, so bitte\\[1.0ex]
\phantom{aaaaaaaaa}ich Sie, mir durch ihn sagen zu lassen, bei wem\\[1.0ex]
\phantom{aaaaaaaaa}ich wohl eine klaine Säge oder ein anderes für\\[1.0ex]
\phantom{aaaaaaaaa}----------------------------------- {\small{new page}} --------------------------------------\\[1.0ex]
\phantom{aaaaaaaaa}meinen Zweck passendes Instrument hier in\\[1.0ex]
\phantom{aaaaaaaaa}Göttingen kaufen lassen kann, und unter\\[1.0ex]
\phantom{aaaaaaaaa}welcher Benennung ich wohl am leichtesten das\\[1.0ex]
\phantom{aaaaaaaaa}Geeignete erhalte.\\[1.0ex]
\phantom{aaaaaaaaa}Genaue Messungen können mir, wie Sie schon\\[1.0ex]
\phantom{aaaaaaaaa}sagten, für den Paukenapparat wenig helfen,\\[1.0ex]
\phantom{aaaaaaaaa}da die Data zu meiner Rechnung sich doch nicht zu{\raisebox{1.5pt}{\tiny=}}\\[1.0ex]
\phantom{aaaaaaaaa}sammenbringen lassen. Aber für die Schnecke\\[1.0ex]
\phantom{aaaaaaaaa}ist eine Rechnung nöthig, und ich hätte gern, wenn\\[1.0ex]
\phantom{aaaaaaaaa}auch nur rohe Angaben über die Abnahme des\\[1.0ex]
\phantom{aaaaaaaaa}Querschnitts der Schneckencanäle und das Breiter\\[1.0ex]
\phantom{aaaaaaaaa}werden der Membrana basilaris von der\\[1.0ex]
\phantom{aaaaaaaaa}Basis gegen die Spitze.\\[1.0ex]
\phantom{aaaaaaaaa}\\[1.0ex]
\phantom{aaaaaaaaaaaaaaaaaaaaaaaaaaaaaa}Ihr\\[1.0ex]
\phantom{aaaaaaaaa}\\[1.0ex]
\phantom{aaaaaaaaa}Göttingen, den 18.\ April\phantom{aaaaaaaaaa}ewig dankbarer\\[1.0ex]
\phantom{aaaaaaaaaaaaaaaa}1866.\phantom{aaaaaaaaaaaaaaaaaaaa}B.\ Riemann\\[1.0ex]
\normalsize
%
%
%
%
%
%
\newpage
\section*{Acknowledgement}
\noindent
I thank Juan M. Marín, Erwin Neuenschwander and Samuel J. Patterson for their support and suggestions. My thanks
also go to Erin C. Rushing (SIL) and Bärbel Mund (SUB) for their friendly permission to use the facsimiles in the
figures. Further, I thank Samuel J. Patterson and my son Wieland for their linguistic corrections.
%
%
%
%
%

%
%
%
%
%
%
\newpage
\section*{Figures}
\label{Figures}
\vspace{8ex}
\piccaption{\label{Abb01}Supplement, first page}
\parpic(\textwidth,128ex){\includegraphics[width=1.00\textwidth]{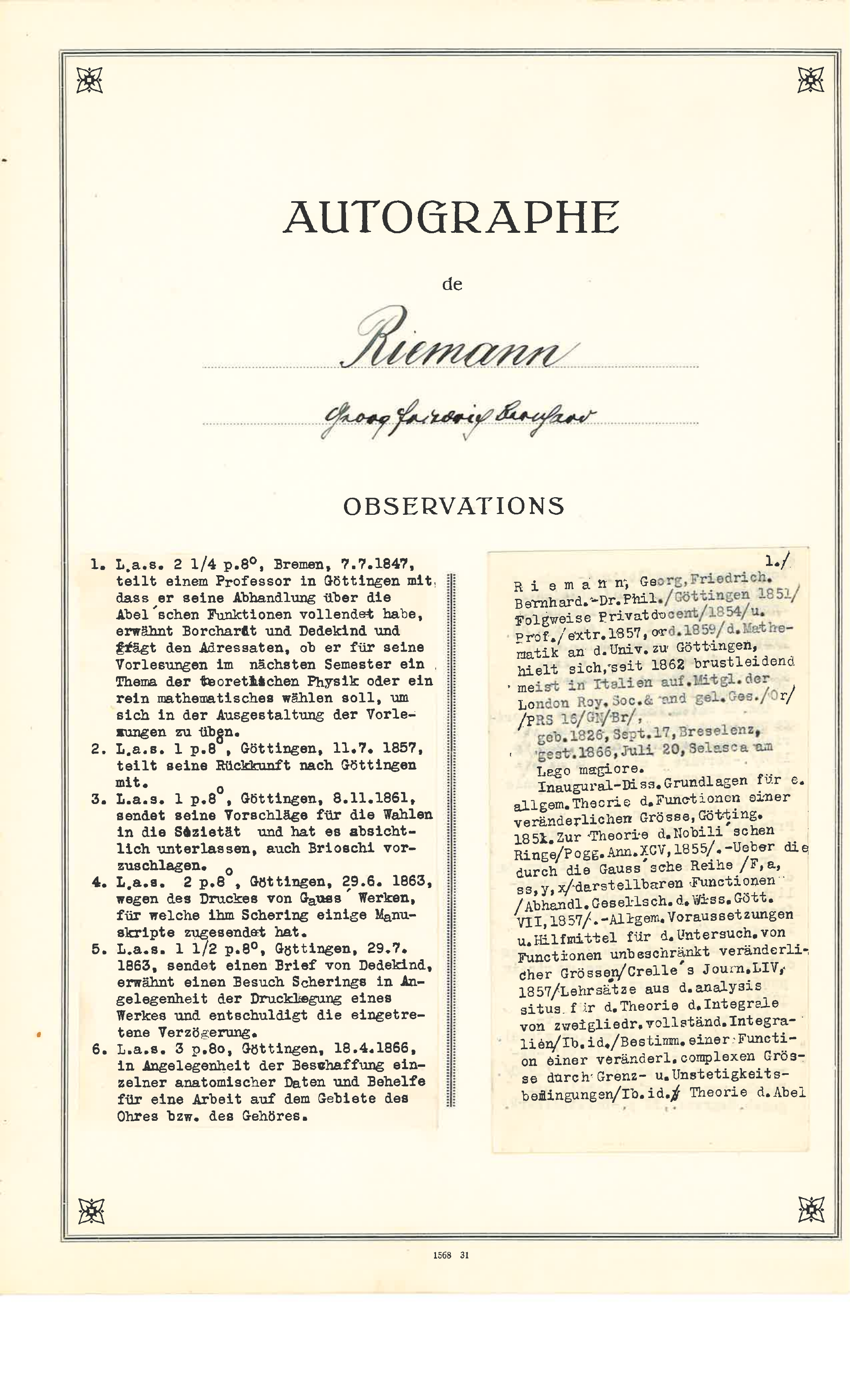}}
\newpage
\phantom{a}
\vspace{8ex}
\piccaption{\label{Abb02}Supplement, second page}
\parpic(\textwidth,129ex){\includegraphics[width=1.00\textwidth]{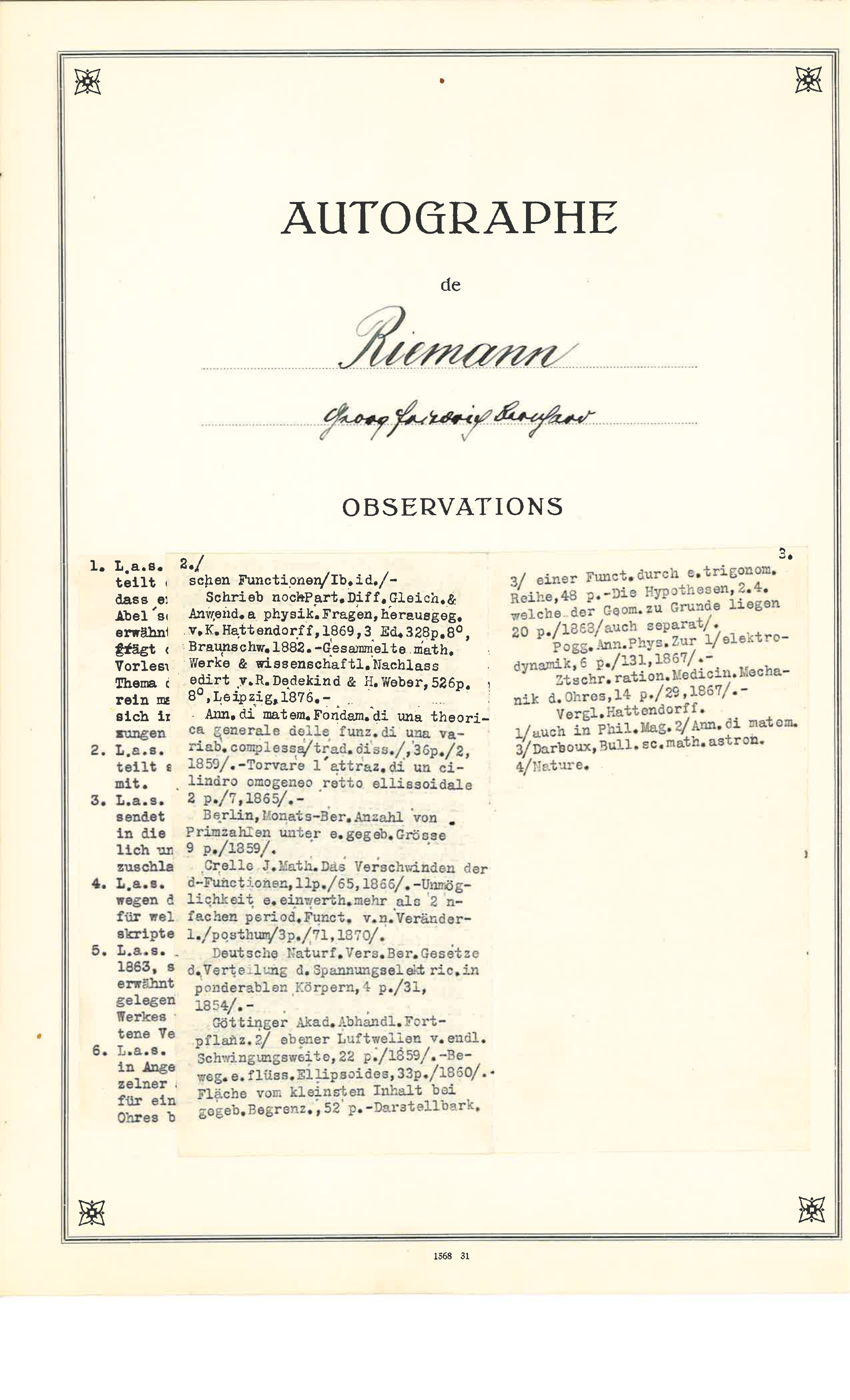}}
\newpage
\piccaption{\label{Abb03}Sixth Letter of 18\textsuperscript{th} April 1866, first page}
\parpic(\textwidth,142ex){\includegraphics[width=1.00\textwidth]{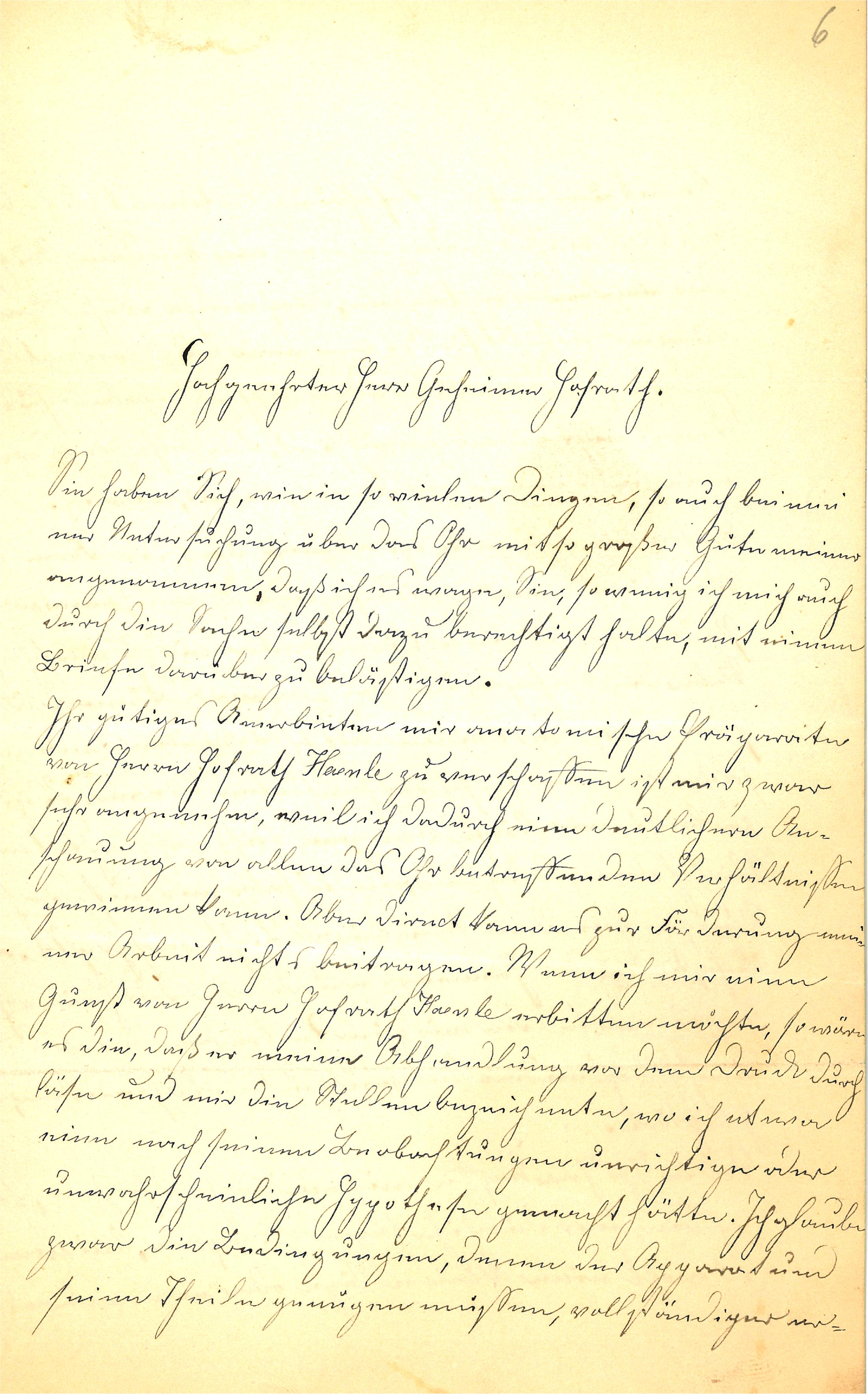}}
\newpage
\piccaption{\label{Abb04}Letter of Elise Riemann to Richard Dedekind, first page}
\parpic(\textwidth,130ex){\includegraphics[width=1.00\textwidth]{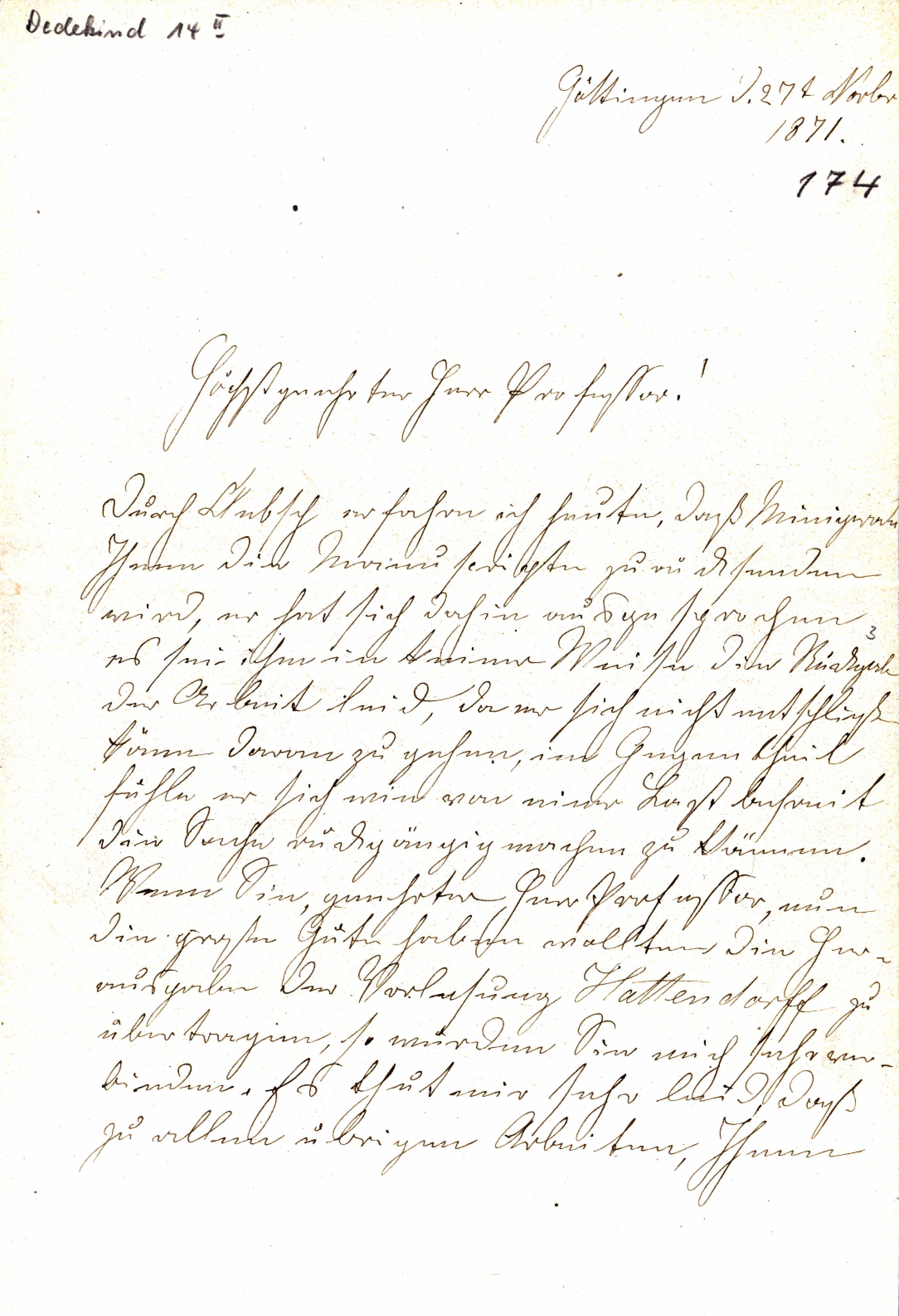}}
\newpage
\piccaption{\label{Abb05}Letter of Elise Riemann to Richard Dedekind, second page}
\parpic(\textwidth,132ex){\includegraphics[width=1.00\textwidth]{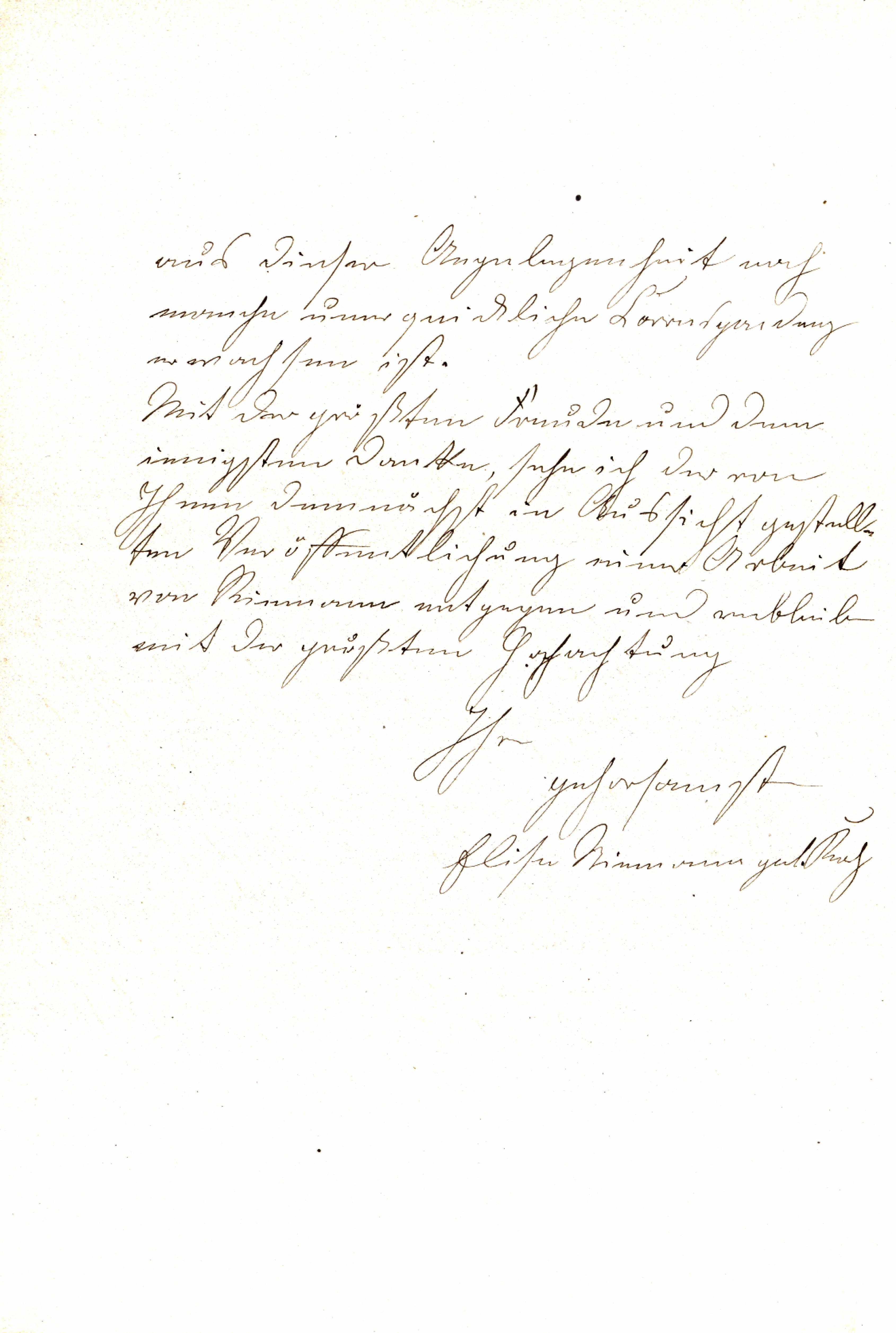}}
%
\end{document}